\documentstyle[epsf,12pt]{article}
\title{The second bounded cohomology \\
of a group with infinitely many ends
}
\author{Koji Fujiwara  \thanks{This work was done when the author visited 
MSRI supported in part by NSF grant DMS-9022140 and a JSPS grant.} \\
Mathematics Department, Keio University \\
Yokohama, 223, Japan \\
{\tt fujiwara@math.keio.ac.jp}}

\date{}
\newtheorem{Th}{Theorem}
\newtheorem{Prop}{Proposition}
\newtheorem{Lemma}{Lemma}[section]
\newtheorem{Cor}{Corollary}

\newcommand{\mathbf}{\bf }
\newcommand{\qed}{\hspace *{\fill} q.e.d.}

\begin{document}
\maketitle

\begin{abstract}
We study the second bounded cohomology of 
an amalgamated free product of groups, and an HNN extension of a group.
As an application, we have a group with 
infinitely many ends has infinite dimensional second bounded cohomology.
\end{abstract}

\noindent
{\bf Contents}

\noindent
\S 1 Introduction

\noindent
\S 2 Quasi homomorphism and the counting function

\noindent 
\S 3 Quasi homomorphisms of $A*_C B$

\noindent
\S 4 Choice of words for $A*_C B$

\noindent
\S 5 Proofs of Th 1 and Cor 1,2,3

\noindent 
\S 6 Quasi homomorphisms of $A*_C \varphi $

\noindent 
\S 7 Choice of words for  $A*_C \varphi $

\noindent 
\S 8 Proofs of Th 2 and Th 3

\section{Introduction}
Bounded cohomology was defined by F. Trauber for groups and by M. Gromov for 
spaces.
We review the definition of bounded cohomology of a group $G$.
Let
$$C^k_b(G;A) = \{ f: G^k \rightarrow A \mid f \mbox{ has bounded
range} \}, $$
where $A={\mathbf Z}$ or ${\mathbf R}$.
The boundary
$\delta:C^k_b(G;A) \rightarrow C^{k+1}_b(G;A)$ is given by
\begin{eqnarray*}
\delta f(g_0,\ldots,g_k) &=& f(g_1,\ldots,g_k)
+ \sum_{i=1}^k (-1)^i f(g_0,\ldots,g_{i-1}g_i,\ldots g_k)\\
 &&\quad + (-1)^{k+1}
f(g_0,\ldots,g_{k-1}).
\end{eqnarray*}
The cohomology of the complex $\{C^k_b(G;A), \delta \}$ is the {\it
bounded cohomology} of $G$, denoted by $H^*_b(G;A)$.
See \cite{G}, \cite{I} as general references for the theory of bounded 
cohomology.

For any group $G$, the first bounded cohomology, $H^1_b(G;A)$
is trivial.
For an amenable group $G$, $H^n_b(G;{\mathbf R})$ is trivial for any $n$.
The first example of a group with non-trivial second bounded cohomology 
was given by R. Brooks, \cite{B}. 
He showed a free group of rank greater than 1 has 
infinite dimensional second bounded cohomology.
R.I. Grigorchuk investigated the structure of the second bounded cohomology
of free groups, torus knot groups and surface groups, \cite{Gr}. 
T. Yoshida \cite{Y} and T. Soma \cite{So1}, \cite{So2} 
studied the third bounded 
cohomology of surfaces and hyperbolic three manifolds.
D.B.A. Epstein and the author showed a non-trivial word-hyperbolic group has 
infinite dimensional second bounded cohomology, \cite{EF}.

We state results on the second bounded cohomology of an amalgamated free
product of groups:

\begin{Th}
Let $G=A*_C B$.
If $|C \backslash A/C| \ge 3$ and $ |B/C| \ge 2$, 
then the cardinality of the dimension of $H^2_b(G;{\mathbf R})$ as 
a vector space over {\bf R} is continuum.
\end{Th}

\begin{Cor}
Let $G=A*B$ with $A \neq \{ 1 \}, B \neq \{ 1 \}$. 
If $G \neq {\mathbf Z}_2 * {\mathbf Z}_2$, 
then the cardinality of the dimension of 
$H^2_b(G;{\mathbf R})$ as a vector space over {\bf R} is continuum.
\end{Cor}

\noindent
{\it Remark}. (1) Cor 1 is a generalization of R. Brooks' result 
on free groups.

\noindent 
(2) Since $\mathbf Z_2 * \mathbf Z_2$ is amenable, 
$H^2_b(\mathbf Z_2 * \mathbf Z_2; \mathbf{R})$ is trivial.

\begin{Cor}
Let $G=A*_C B$. If $|A|=\infty$, $|C| < \infty$, and $ |B/C| \ge 2$,
then the cardinality of the dimension of 
$H^2_b(G;\mathbf R)$ as a vector space over {\bf R} is continuum.
\end{Cor}

\begin{Cor}
Let $G=A*_C B$. If $A$ is abelian, $|A/C| \ge 3$, and $ |B/C| \ge 2$,
then the cardinality of the dimension of 
$H^2_b(G;\mathbf R)$ as a vector space over {\bf R} is continuum.
\end{Cor}

\noindent
{\it Example}. $PSL_2(\mathbf Z)=\mathbf Z_2 * \mathbf Z_3$
and  $SL_2(\mathbf Z)=\mathbf Z_4*_{\mathbf Z_2} \mathbf Z_6$ 
satisfy the assumption of Th 1. They are non-elementary word-hyperbolic 
groups too.

For an HNN extension of a group, we show:

\begin{Th}
Let $G=A*_C \varphi$. If $|A/C| \ge 2, \, |A/ \varphi (C)| \ge 2$, then 
then the cardinality of the dimension of 
$H^2_b(G;\mathbf R)$ as a vector space over {\bf R} is continuum.
\end{Th}

\noindent
{\it Remark}. In Th 1 and 2, and Cor 1, 2 and 3, we don't have to assume 
$G$ is finitely generated.

We apply Th 1 and 2 to a group with infinitely many ends, namely, 
due to Stallings' structure theorem \cite{S}, we have:

\begin{Th}
If $G$ is a finitely generated group with infinitely many ends, 
then the cardinality of the dimension of $H^2_b(G;\mathbf R)$ as a vector 
space over {\bf R} is continuum.
\end{Th}

To conclude, we state a conjecture.

\noindent
{\it Conjecture}. 
Let $G$ be a group. If $H_b^2(G;{\mathbf R}) \neq \{0\}$, 
then $H_b^2(G;{\mathbf R})$ is infinite dimensional as a vector space over 
{\bf R}.

\section{Quasi homomorphism and the counting function}

We will review the counting function of a group w.r.t. a word. 
This was 
defined in \cite{EF} as a generalization of R. Brooks' counting function 
for free groups in \cite{B}. 
Let $G$ be a group with a (finite or infinite) set of generators and 
$\Gamma(G)$ the Cayley graph.
For a word $w=x_1 x_2 \dots x_n$ in these generators, 
define $|w|=n$. 
Let $\overline{w}$ be the element of $G$ which is represented 
by the word $w$. 
Define $w^{-1} =x_n^{-1} \dots x_1^{-1}.$
We sometimes identify a word $w$ and the path starting 
at $1$ and labeled by $w$ in $\Gamma(G)$. 
For a path $\alpha$ labeled by $w$, 
define $|\alpha|=|w|$ and $\overline{\alpha}=\overline{w}$.
For an element $g$ in $G$, 
define $|g|=\inf_{\alpha} |\alpha|$, 
where $\alpha$ ranges over all the paths 
with $\overline \alpha =g$.

Let $\alpha$ be a finite path in $\Gamma(G)$. Define 
$|\alpha|_w$ to be the maximal number of times that $w$ can 
be seen as a subword 
of $\alpha$ without overlapping. We define 

$$c_w(\alpha) = \sup _{\alpha'} \{ |\alpha'|_w 
                        - ( |\alpha'|-|\overline{\alpha}|)\}
              = |\overline{\alpha}| - \inf_{\alpha'}(|\alpha'|-|\alpha'|_w), $$
where $\alpha'$ ranges over all the paths with the same starting point 
as $\alpha$ and the same finishing point.
If the infimum in the definition of $c_w(\alpha)$ is attained by 
$\alpha'$, we say that $\alpha'$ {\it realizes} $c_w$ {\it at} $\alpha$.
We have $|\overline{\alpha} | / |w| \ge c_w(\alpha) \ge 0$.
If $\alpha$ is a geodesic, then $c_w(\alpha) \ge |\alpha|_w$.

We define $h_w =c_w - c_{w^{-1}}$. 
For each $g$ in $G$, we choose $\gamma_g$ to be a path from 1
to $g$ and set $h_w(g)=h_w(\gamma_g)$. 
Then $h_w(g)$ does not depend on the choice of $\gamma_g$. 
Thus $h_w \in C^1(G;\mathbf Z)$.

Let $f \in C^1(G;\mathbf R)$. If there exists a constant $D < \infty$ s.t.
$$|f(gh) -f(g) -f(h) | \le D,$$
for any $g,h \in G$, then we say $f$ is a {\it quasi homomorphism} with 
{\it defect} $D$. 
Let $f$ be a quasi homomorphism with defect $D$. 
Then $|\delta f| \le D$ and $\delta (\delta f)=0$, thus 
$\delta f \in Z^2_b(G;\mathbf R)$, which defines 
$[\delta f] \in H^2_b(G;\mathbf R)$.
Remark that we always have $[\delta f] =0$ in $H^2(G; \mathbf R)$, however,
we may have $[\delta f] \not = 0$ in $H^2_b(G;\mathbf R)$ since 
$f$ is not necessarily in $C^1_b(G;\mathbf R)$.

\section{Quasi homomorphisms of $A*_C B$}

Let $G=A*_C B$ with $ |A/C| \ge 2, \, |B/C| \ge 2$.
Take the set $\{A \cup B \} \backslash \{1 \}$ as a set of 
generators of $G$ and denote its Cayley graph by $\Gamma(G)$.
Note that if a generating set is infinite, then $\Gamma$ is 
not locally finite.

If a word $w= x_1 \dots x_n$ satisfies $n=1$ or 
$x_1, x_3, \dots \in A \backslash C$
(or $B \backslash C$) and $x_2, x_4, \dots \in B \backslash C$
(or $A \backslash C$, resp.), then we say $w$ is {\it reduced}.

\begin{Lemma}
A word $w =x_1 \dots x_n$ is reduced iff it is a geodesic in $\Gamma$.
\end{Lemma}

\noindent
{\it Proof}. Assume $w$ is not reduced, then there exists a subword 
$x_ix_{i+1}$ with $x_i,x_{i+1} \in A$(or $x_i,x_{i+1} \in B$). 
Then $\overline{x_ix_{i+1}} \in A(\mbox{or }B \mbox{ resp.})$, thus 
$|\overline{x_ix_{i+1}}| \le 1$.
Therefore $x_ix_{i+1}$ is 
not a geodesic, hence $w$ is not a geodesic. 

On the other hand, assume $w$ is not a geodesic.
To show $w$ is not reduced by contradiction, suppose $w$
is reduced.
Take a geodesic $\gamma$, hence reduced,  s.t. 
$\overline{w} = \overline{\gamma}$. 
Then $|w| = |\gamma|$ since reduced words representing a same 
element have same length. Thus $w$ is a geodesic.
This is a contradiction. We showed $w$ is not reduced.
\qed

\begin{Lemma}
Let $w$ be a word and $\alpha$ a path. If $w^2$ is reduced, 
then there is a geodesic which realizes $c_w$ at $\alpha$.
\end{Lemma}

\noindent
{\it Proof}. Since $w^2$ is reduced, $w$ is reduced. 
Let $\gamma$ be a path which realizes $c_w$ at $\alpha$ 
s.t. $|\gamma|_w$ is minimal among all the realizing paths at $\alpha$.
We will show that $\gamma$ is a geodesic.
If $|\gamma|_w =0$, then $\gamma$ is a geodesic.
Assume $|\gamma|_w =n>0$. 
Write $\gamma$ as 
$$\gamma_1 w_1 \gamma_2 \dots w_n \gamma_{n+1},$$
where $w_i$ is a copy of $w$ and $\gamma_i$ may be an empty word.
First, to show every $\gamma_i$ is reduced by contradiction, 
suppose some $\gamma_i$ is not reduced.
Replace $\gamma_i$ by 
a reduced word $\gamma'_i$ with $\overline{\gamma_i}=\overline{\gamma'_i}$,
then we have a new path
$$\gamma' = \gamma_1 w_1 \gamma_2 \dots \gamma'_i \dots w_n \gamma_{n+1},$$
which satisfies $|\gamma'| < |\gamma|$, $|\gamma'|_w \ge |\gamma|_w$
and $\overline{\gamma}=\overline{\gamma'}$. 
This contradicts that $\gamma$ is a realizing path.
Thus every $\gamma_i$ is reduced.
Next, in order to show that $\gamma$ is reduced by contradiction, suppose not.
Since $w^2$ is reduced, there is a subword 
$w_i \gamma _{i+1} w_{i+1}$ of $\gamma$ which is not reduced and 
$\gamma_{i+1}$ is not empty. 
Since $w_i, w_{i+1}, \gamma_{i+1}$ are reduced, one of the following 
four cases is occurring.

\noindent
(i) The last letter of $w_i$ and the initial letter of $\gamma_{i+1}$
     are in $A$.

\noindent
(ii) The last letter of $w_i$ and the initial letter of $\gamma_{i+1}$
     are in $B$.

\noindent
(iii) The last letter of $\gamma_{i+1}$ and the initial letter 
       of $w_{i+1}$ are in $A$.

\noindent
(iv) The last letter of $\gamma_{i+1}$ and the initial letter 
       of $w_{i+1}$ are in $B$.

We treat the case (i).
Write 
$$w_i = \dots b_1 a_1, \, \gamma_{i+1} = a_2 b_2 \dots, \, a_i \in A,
 b_i \in B.$$
Rewrite the subword $w_i \gamma_{i+1}$ in $\gamma$ as 
$$w_i \gamma_{i+1} = \dots b_1 a_1 a_2 b_2 \dots = \dots b_1 a' b_2 \dots,$$
where $a' =a_1a_2 \in A.$
This gives a new word $\gamma'$ with $|\gamma'|=|\gamma|-1,$
$|\gamma'|_w = |\gamma|_w -1$. 
Since $\overline{\gamma'}= \overline{\gamma}$ and 
$|\gamma|- |\gamma|_w = |\gamma'|-|\gamma'|_w ,$ $\gamma'$ is another 
realizing path with $|\gamma'|_w < |\gamma|_w $. 
This contradicts the choice of $\gamma$, thus $\gamma$ is reduced. 
By Lemma 3.1, $\gamma$ is a geodesic.
A similar argument applies to the other three cases, and we omit it.
\qed

\begin{Lemma}
Let $\alpha, \beta$ be paths starting at $1$. We have 
$$|c_w(\alpha)-c_w(\beta)| \le 2|\overline{\alpha^{-1} \beta}|.$$
\end{Lemma}

\begin{figure}
\centerline{
\epsfbox{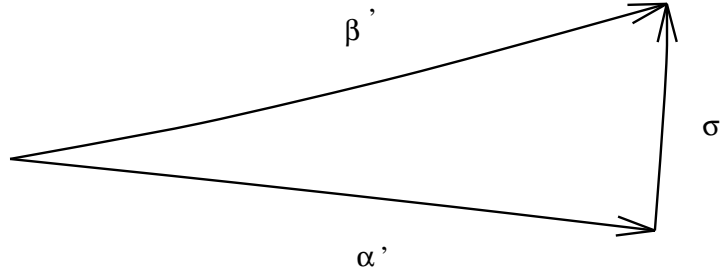}
}
\caption{Geodesic triangle. This illustrates Lemma 3.3.}
\end{figure}

\noindent
{\it Proof}. Take realizing geodesics of $c_w$ at $\alpha$, $\alpha'$, and 
at $\beta$, $\beta'$. Then $c_w(\alpha)=|\alpha'|_w, \, 
c_w(\beta)=|\beta'|_w$. 
Take a geodesic $\sigma$ with 
$\overline{\alpha^{-1} \beta} =\overline{\sigma}$.

Since the path $\alpha' \sigma$ satisfies 
$\overline{\alpha' \sigma}=\overline{\beta'}$, 
$$|\beta'|-|\beta'|_w \le |\alpha' \sigma|-|\alpha' \sigma|_w.$$
Since $|\alpha' \sigma| = |\alpha'| +|\sigma|$ and
$|\alpha' \sigma|_w \ge |\alpha'|_w + |\sigma|_w,$
$$|\beta'|-|\beta'|_w \le |\alpha'|+ |\sigma| -|\alpha'|_w - |\sigma|_w. $$
Thus 
$$
c_w(\beta)=|\beta'|_w \ge |\alpha'|_w + |\beta'| -|\alpha'| -|\sigma|
                             +|\sigma|_w 
                     \ge c_w(\alpha) -2|\sigma|,$$
since $|\alpha'|_w =c_w(\alpha), \, |\beta'| -|\alpha'| \ge - |\sigma|, \,
       |\sigma|_w \ge 0.$
Similarly, $c_w(\alpha) \ge c_w(\beta) - 2|\sigma|.$
Thus $|c_w(\alpha) -  c_w(\beta)| \le 2|\sigma| 
=2|\overline{\alpha^{-1} \beta}|.$
\qed

\begin{Lemma}
Let $\alpha, \beta$ be paths starting at $1$. We have 
$$|h_w(\alpha)-h_w(\beta)| \le 4| \overline{\alpha^{-1} \beta} | .$$
\end{Lemma}

\noindent
{\it Proof}. $h_w =c_w -c_{w^{-1}}$ and Lemma 3.3. \qed

\begin{Lemma}
$$c_w(\alpha) =c_{w^{-1}}(\alpha^{-1}).$$
\end{Lemma}

\noindent
{\it Proof}. Clear from the definition of $c_w$. \qed

\begin{Lemma}
$$h_w(\alpha) = -h_w(\alpha^{-1}).$$
\end{Lemma}

\noindent
{\it Proof}.
$h_w(\alpha^{-1})=c_w(\alpha^{-1})-c_{w^{-1}}(\alpha^{-1})
=c_{w^{-1}}(\alpha) - c_w(\alpha) = -h_w(\alpha)$. \qed

\begin{Lemma}
Let $\alpha$ be a reduced path. If $\alpha=\alpha_1 \alpha_2$, then 
$$|h_w(\alpha) - h_w(\alpha_1) -  h_w(\alpha_2)| \le 10.$$
\end{Lemma}

\begin{figure}
\centerline{
\epsfbox{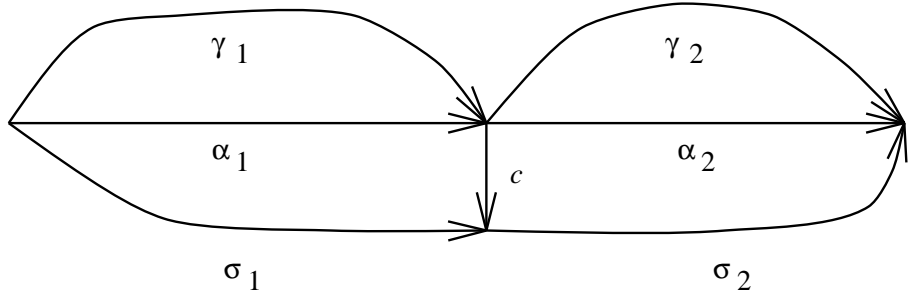}
}
\caption{Dividing geodesics. This illustrates Lemma 3.7. 
$\alpha =\alpha_1 \alpha_2$, 
$\gamma = \gamma_1 \gamma_2$, $\sigma=\sigma_1 \sigma_2$.}
\end{figure}

\noindent
{\it Proof}. We will show 
$$|c_w(\alpha) - c_w(\alpha_1) -  c_w(\alpha_2)| \le 5.$$
Take realizing geodesics of $c_w$ at $\alpha_1$, $\gamma_1$, and 
at $\alpha_2$, $\gamma_2$.
We have $|\gamma_1|_w =c_w(\alpha_1), \, |\gamma_2|_w =c_w(\alpha_2).$
Put $\gamma = \gamma_1 \gamma_2$, then $\overline{\gamma} =\overline{\alpha}$.
Since $\alpha$ is reduced, it is a geodesic. 
Since $|\alpha|=|\alpha_1|+|\alpha_2|=|\gamma_1|+|\gamma_2|=|\gamma|$, 
$\gamma$ is a geodesic.
Since $\gamma$ is a geodesic with $\overline{\gamma}=\overline{\alpha}$,
$c_w(\alpha) \ge |\gamma|_w$. Thus 
$$c_w(\alpha) \ge |\gamma|_w \ge |\gamma_1|_w + |\gamma_2|_w = 
c_w(\alpha_1) + c_w(\alpha_2).$$

On the other hand, 
take a realizing geodesic $\sigma$ at $\alpha$, then $c_w(\alpha)=|\sigma|_w$. 
Since $\alpha$ and $\sigma$ are reduced, there exists subdivision of $\sigma$,
$\sigma =\sigma_1 \sigma_2$ s.t. 
$\overline{\alpha_1^{-1} \sigma_1} = \overline{\alpha_2 \sigma_2^{-1}} =c, $
for some $c \in C$.
Since $|c| \le 1$, from Lemma 3.3 and 3.5.
$$|c_w(\alpha_i)-c_w(\sigma_i)| \le 2, \, i=1,2.$$
Since $\sigma_1, \sigma_2$ are geodesics, 
$|\sigma_i|_w \le c_w(\sigma_i), \, i=1,2.$
Thus
$$c_w(\alpha)= |\sigma|_w \le |\sigma_1|_w + |\sigma_2|_w +1 \le 
c_w(\sigma_1) + c_w(\sigma_2) +1 \le c_w(\alpha_1) + c_w(\alpha_2) +5.$$
We showed $|c_w(\alpha) - c_w(\alpha_1) -c_w(\alpha_2)| \le 5$.
$|c_{w^{-1}}(\alpha) - c_{w^{-1}}(\alpha_1) -c_{w^{-1}}(\alpha_2)| \le 5$ 
is similar. Since $h_w=c_w -c_{w^{-1}}$, we get 
$|h_w(\alpha) - h_w(\alpha_1) -  h_w(\alpha_2)| \le 10.$ 
\qed

\begin{Prop}
Let $w$ be a word. We have  
$$|\delta h_w| \le 78. $$
\end{Prop}

\begin{figure}
\centerline{
\epsfbox{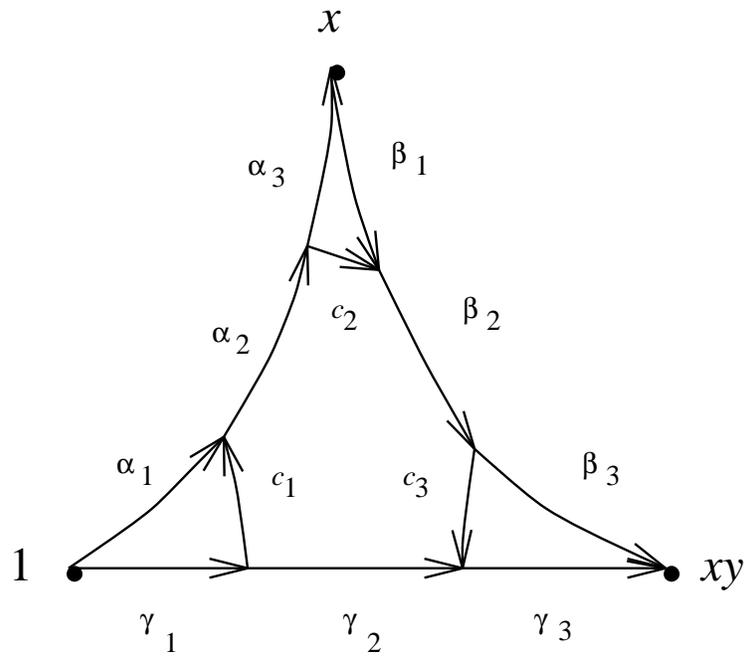}
}
\caption{Dividing geodesic triangle. This illustrates Prop 1. 
$c_1, c_2, c_3 \in C$.}
\end{figure}

\noindent
{\it Proof}.
Let $x,y$ be elements in $G$. 
$\delta h_w(x,y) =h_w(x) + h_w(y) -h_w(xy).$
We will show 
$|h_w(xy)-h_w(x)-h_w(y)| \le 78.$
Take geodesics, hence reduced paths, $\alpha, \beta$ and $\gamma$ with 
$\overline{\alpha} =x, \overline{\beta} =y, \overline{\gamma} =xy$.
Since $\overline{\alpha \beta} = \overline{\gamma}$ and 
$\alpha, \beta, \gamma$ are 
reduced, there exist subdivisions of $\alpha, \beta, \gamma$,
$$\alpha= \alpha_1 \alpha_2 \alpha_3, \, 
  \beta =\beta_1 \beta_2 \beta_3, \, 
  \gamma =\gamma_1 \gamma_2 \gamma_3, $$
s.t. 
$$\overline{\gamma_1^{-1} \alpha_1} =c_1, \, 
  \overline{\alpha_3 \beta_1} =c_2, \, 
  \overline{\beta_3 \gamma_3^{-1}} =c_3, $$
for some $c_1, c_2, c_3 \in C$ and $\overline{\alpha_2}, \overline{\beta_2}, 
\overline{\gamma_2}$ are (simultaneously) in $A$ or $B$.

Since $\overline{\alpha_2}$ is in $A$ or $B$, 
$c_w(\alpha_2), c_{w^{-1}}(\alpha_2) \le |\overline{\alpha_2}| \le 1$. 
Thus $|h_w(\alpha_2)| \le 2$. Similarly, 
$|h_w(\beta_2)|, |h_w(\gamma_2)| \le 2$.

Using Lemma 3.7 twice for $\alpha = ((\alpha_1 \alpha_2) \alpha_3)$, 
$$|h_w(\alpha) -h_w(\alpha_1) -h_w(\alpha_2)-h_w(\alpha_3)| \le 20, $$
thus $|h_w(\alpha) -h_w(\alpha_1) -h_w(\alpha_3) | \le 22.$
Similarly, 
$$|h_w(\beta) - h_w(\beta_1) - h_w(\beta_3)| \le 22, $$
$$|h_w(\gamma) - h_w(\gamma_1) - h_w(\gamma_3)| \le 22. $$
Thus
\begin{eqnarray*}
&&         |h_w(\alpha) + h_w(\beta) -h_w(\gamma)| \\
& \le &   |h_w(\alpha_1) +h_w(\alpha_3) + h_w(\beta_1) + h_w(\beta_3) 
          -h_w(\gamma_1) - h_w(\gamma_3) | + 66\\
& \le &    |h_w(\alpha_1) - h_w(\gamma_1)| 
          + | h_w(\beta_1) - h_w(\alpha_3^{-1})|
          + |h_w(\gamma_3^{-1}) - h_w(\beta_3^{-1})| + 66 \\
&\le &    4(|c_1| +|c_2| +|c_3|) + 66 \le 78,
\end{eqnarray*}
since $|c_i| \le 1$ and Lemma 3.4 and 3.6.

By definition, $h_w(\alpha)=h_w(x), h_w(\beta)=h_w(y), h_w(\gamma)=h_w(xy)$, 
thus 
$$|h_w(x) + h_w(y) -h_w(xy) | \le 78.$$
\qed

\section{Choice of words for $A*_C B$}

The goal of this section is to prove the following proposition.
\begin{Prop}
Let $G=A*_CB$ with $|C \backslash A /C| \ge 3$ and $|B/C| \ge 2$.
Then there exist words $w_i, \, 0 \le i < \infty$ which satisfy
the following conditions.

\noindent
(1) For any $i \ge 0$ and any $n \ge 1$, $h_{w_i}(w_i^n) =n$.

\noindent
(2) For any $j > i \ge 0$ and any $n \ge 1$, 
$h_{w_j}(w_i^n) =0$.

\noindent
(3) For any $i \ge 0$, 
$\overline{w_i} \in [G,G]$.

\end{Prop}

Let $G=A*_CB$ with $|C \backslash A /C| \ge 3$ and $|B/C| \ge 2$.
Let $w$ and $w'$ be reduced paths starting at a same point and 
finishing at a same point,
$$w=a_1 b_1 \dots a_n b_n, $$
$$w'=a_1' b_1' \dots a_m' b_m',$$
where $a_1, b_n, a_1', b_m'$ may be empty.

Since $\overline{w}=\overline{w'}$ and $w$ and $w'$ are reduced,
$n=m$ and there exist $c_1, \dots c_n$ and $c_0', c'_1, \dots c'_n$ in $C$
s.t. 
$$c_{i-1}' a_i c_i^{-1} = a_i', \, c_ib_i c_i'^{-1} = b_i',$$
for $1\le ^{\forall}i \le n$, where $c_0'=c_n'=1$. See Figure 4.

\begin{figure}
\centerline{\epsfxsize=12cm
\epsfbox{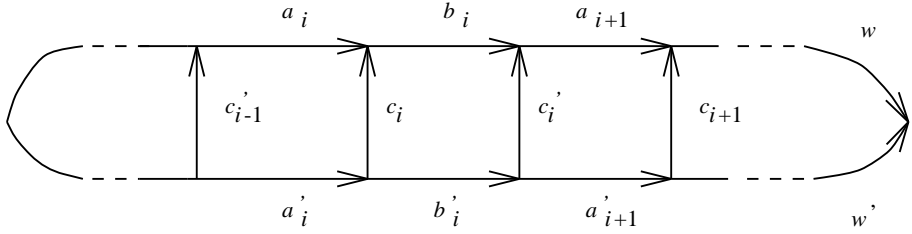}
}
\caption{Reduced paths $w, w'$ with common end points.
$\overline{w} = \overline{w'}$.}
\end{figure}

Let $v$ be a subword of $w$,
$$v=a_ib_i \dots a_j b_j,$$
where $a_i, b_j$ may be empty. 
Define a subword of $w'$, denoted by $P(v)$, by 
$$P(v) =a_i' b_i' \dots a_j' b_j'.$$

\begin{figure}
\centerline{\epsfxsize=12cm
\epsfbox{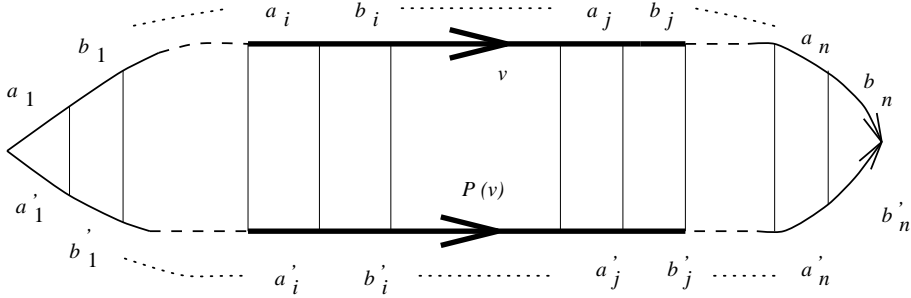}
}
\caption{Definition of $P(v)$.}
\end{figure}

Let $v'$ be a subword of $w'$.
If $v'$ is a subword of $P(v)$, we say 
$v$ {\it covers} $v'$. 
If $v'=P(v)$, then we say $v$ {\it faces} $v'$.

Since $|C \backslash A / C| \ge 3$, we can choose $a_1, a_2 \in A$ s.t.
$a_1, a_2 \in A \backslash C,$ and $a_2 \not \in Ca_1C$.  
Take $ b \in B \backslash C$. 
Define words $w_i, 0 \le i < \infty$ by 
\begin{eqnarray}
w_i &=& (a_1b)^{10^i} (a_1^{-1}b^{-1})^{10^i} 
      (a_2b)^{10^i} (a_2^{-1}b^{-1})^{10^i}    \label{words} \nonumber \\ 
    && (a_1b)^{4 \cdot 10^i} (a_1^{-1}b^{-1})^{4 \cdot 10^i}
      (a_2b)^{4 \cdot 10^i} (a_2^{-1}b^{-1})^{4 \cdot 10^i}. \nonumber
\end{eqnarray}
We write the subword $(a_1b)^{4 \cdot 10^i}$ by $w_i(1,+)$ and 
the subword $(a_1^{-1}b^{-1})^{4 \cdot 10^i}$ by $w_i(1,-)$. 

\begin{Lemma}
The words $w_i, 0 \le i < \infty$ satisfy the following conditions.

\noindent
(1) For any $i \ge 0$ and any $n \ge 1$, 
$w_i^n$ is reduced.

\noindent
(2) For any $i \ge 0$, 
$|w_i|=40 \cdot 10^i$ and $|w_i(1,\pm)| =8 \cdot 10^i.$

\noindent
(3) For any $i \ge 0$, 
$\overline{w_i} \in [G,G].$
\end{Lemma}

\noindent
{\it Proof}. (1) is clear from the choice of $a_1, a_2$, and $b$.
(2) is obvious. For (3), 
use $ \overline{(a_ib)^p(a_i^{-1}b^{-1})^p} \in [G,G]$ 
for $i=1,2$ and any $p \ge 1$. \qed

\begin{Lemma}
In any pair of paths having a same starting point and a same finishing point, 
the following conditions hold.

\noindent
(1) $a_1$ cannot face $a_2$. 

\noindent
(2) For any $i \ge 0$, 
$w_i^2$ cannot cover $w_i^{-1}$.

\noindent
(3) For any $i < j$, $w_j(1,+)$ cannot cover $w_i$ and 
$w_j(1,-)^{-1}$ cannot cover $w_i$.

\noindent
(4) For any $k>0$ and any $i<j$, 
$w_i^k$ cannot cover $w_j$ nor $w_j^{-1}$.
\end{Lemma}

\noindent
{\it Proof}. 
(1) If $a_1$ faces $a_2$ in some pair of paths, 
then there exist $c_1, c_2 \in C$ s.t.
$c_1a_1c_2 =a_2$, thus $a_2 \in Ca_1C$, which contradicts our choice 
of $a_1$ and $a_2$.

\noindent
(2) Let $W_0$ symbolize the order of the occurrence of elements 
labeled by $a_1^{\pm 1}$ or  $a_2^{\pm 1}$ in $w_0$ in the following way:
$$W_0=1 \overline{1} 2 \overline{2}
         1111  \overline{1}  \overline{1}  \overline{1}  \overline{1} 
         2222  \overline{2} \overline{2} \overline{2}  \overline{2}, $$
where $1$ stands for $a_1$, 2 for $a_2$, $\overline{1}$ for $a_1^{-1}$,
and $\overline{2}$ for $a_2^{-1}$. 
Here, remark that in $W_0$, we ignore $b^{\pm 1}$'s of $w_0$.
Let $W_0^{-1}$ represent for $w_0^{-1}$:
$$W_0^{-1}= 2222 \overline{2} \overline{2} \overline{2}  \overline{2}
            1111 \overline{1}  \overline{1}  \overline{1}  \overline{1} 
            2  \overline{2} 1 \overline{1}   .$$

In order to show the conclusion by contradiction, 
assume $w_0^2$ covers $w_0^{-1}$. 
Then $W_0^2$ covers $W_0^{-1}$. 
For each possible position for $W_0^2$ covering $W_0^{-1}$, one can find 
some $1$(or $2$, $\overline{1}$, $\overline{2}$) in $W_0^{-1}$ 
facing some $2$ (or $1$, $\overline{2}$, $\overline{1}$ resp.) in $W_0^2$,
see Figure 6.
This means $a_1$ is facing some $a_2$, which contradicts (1).
Thus we showed (2)
for $i=0$. A similar argument works for $w_i$ for any $i \ge 1$, 
and we omit it.

\begin{figure}
\begin{center}
\setlength{\unitlength}{0.6mm}
\begin{picture}(220,220)(-10,-200)
\put(-15,0){$W_0^2$:}
\put(0,0){1}
\put(5,0){$\overline{1}$}
\put(10,0){2}
\put(15,0){$\bar{2}$}
\put(20,0){1}
\put(25,0){1}
\put(30,0){1}
\put(35,0){1}
\put(40,0){$\overline{1}$}
\put(45,0){$\overline{1}$}
\put(50,0){$\overline{1}$}
\put(55,0){$\overline{1}$}
\put(60,0){2}
\put(65,0){2}
\put(70,0){2}
\put(75,0){2}
\put(80,0){$\overline{2}$}
\put(85,0){$\overline{2}$}
\put(90,0){$\overline{2}$}
\put(95,0){$\overline{2}$}

\put(100,0){1}
\put(105,0){$\overline{1}$}
\put(110,0){2}
\put(115,0){$\overline{2}$}
\put(120,0){1}
\put(125,0){1}
\put(130,0){1}
\put(135,0){1}
\put(140,0){$\overline{1}$}
\put(145,0){$\overline{1}$}
\put(150,0){$\overline{1}$}
\put(155,0){$\overline{1}$}
\put(160,0){2}
\put(165,0){2}
\put(170,0){2}
\put(175,0){2}
\put(180,0){$\overline{2}$}
\put(185,0){$\overline{2}$}
\put(190,0){$\overline{2}$}
\put(195,0){$\overline{2}$}

\put(0,-10){2}  \put(0,-14){$\bullet$}              
\put(5,-10){2}
\put(10,-10){2}
\put(15,-10){2}
\put(20,-10){$\overline{2}$}
\put(25,-10){$\overline{2}$}
\put(30,-10){$\overline{2}$}
\put(35,-10){$\overline{2}$}
\put(40,-10){1}
\put(45,-10){1}
\put(50,-10){1}
\put(55,-10){1}
\put(60,-10){$\overline{1}$}
\put(65,-10){$\overline{1}$}
\put(70,-10){$\overline{1}$}
\put(75,-10){$\overline{1}$}
\put(80,-10){2}
\put(85,-10){$\overline{2}$}
\put(90,-10){1}
\put(95,-10){$\overline{1}$} \put(95,-14){$\bullet$}  
\put(5,-20){2}
\put(10,-20){2}
\put(15,-20){2}
\put(20,-20){2} \put(20,-24){$\bullet$}
\put(25,-20){$\overline{2}$}
\put(30,-20){$\overline{2}$}
\put(35,-20){$\overline{2}$}
\put(40,-20){$\overline{2}$}  \put(40,-24){$\bullet$}
\put(45,-20){1}
\put(50,-20){1}
\put(55,-20){1}
\put(60,-20){1} \put(60,-24){$\bullet$}
\put(65,-20){$\overline{1}$}
\put(70,-20){$\overline{1}$}
\put(75,-20){$\overline{1}$}
\put(80,-20){$\overline{1}$} \put(80,-24){$\bullet$}
\put(85,-20){2}
\put(90,-20){$\overline{2}$}
\put(95,-20){1}
\put(100,-20){$\overline{1}$}
\put(10,-30){2}
\put(15,-30){2}
\put(20,-30){2} \put(20,-34){$\bullet$}
\put(25,-30){2} \put(25,-34){$\bullet$}
\put(30,-30){$\overline{2}$}
\put(35,-30){$\overline{2}$}
\put(40,-30){$\overline{2}$}  \put(40,-34){$\bullet$}
\put(45,-30){$\overline{2}$}  \put(45,-34){$\bullet$}
\put(50,-30){1}
\put(55,-30){1}
\put(60,-30){1}  \put(60,-34){$\bullet$}
\put(65,-30){1}  \put(65,-34){$\bullet$}
\put(70,-30){$\overline{1}$}
\put(75,-30){$\overline{1}$}
\put(80,-30){$\overline{1}$}  \put(80,-34){$\bullet$}
\put(85,-30){$\overline{1}$}  \put(85,-34){$\bullet$}
\put(90,-30){2}
\put(95,-30){$\overline{2}$}
\put(100,-30){1}
\put(105,-30){$\overline{1}$}
\put(15,-40){2}
\put(20,-40){2} \put(20,-44){$\bullet$}
\put(25,-40){2} \put(25,-44){$\bullet$}
\put(30,-40){2} \put(30,-44){$\bullet$}
\put(35,-40){$\overline{2}$}
\put(40,-40){$\overline{2}$}  \put(40,-44){$\bullet$}
\put(45,-40){$\overline{2}$}  \put(45,-44){$\bullet$}
\put(50,-40){$\overline{2}$}  \put(50,-44){$\bullet$}
\put(55,-40){1}
\put(60,-40){1}  \put(60,-44){$\bullet$}
\put(65,-40){1}  \put(65,-44){$\bullet$}
\put(70,-40){1}  \put(70,-44){$\bullet$}
\put(75,-40){$\overline{1}$}
\put(80,-40){$\overline{1}$}  \put(80,-44){$\bullet$}
\put(85,-40){$\overline{1}$}  \put(85,-44){$\bullet$}
\put(90,-40){$\overline{1}$}  \put(90,-44){$\bullet$}
\put(95,-40){2}
\put(100,-40){$\overline{2}$}
\put(105,-40){1}
\put(110,-40){$\overline{1}$}
\put(20,-50){2} \put(20,-54){$\bullet$}
\put(25,-50){2} \put(25,-54){$\bullet$}
\put(30,-50){2} \put(30,-54){$\bullet$}
\put(35,-50){2} \put(35,-54){$\bullet$}
\put(40,-50){$\overline{2}$} \put(40,-54){$\bullet$}
\put(45,-50){$\overline{2}$} \put(45,-54){$\bullet$}
\put(50,-50){$\overline{2}$} \put(50,-54){$\bullet$}
\put(55,-50){$\overline{2}$} \put(55,-54){$\bullet$}
\put(60,-50){1} \put(60,-54){$\bullet$}
\put(65,-50){1} \put(65,-54){$\bullet$}
\put(70,-50){1} \put(70,-54){$\bullet$}
\put(75,-50){1} \put(75,-54){$\bullet$}
\put(80,-50){$\overline{1}$} \put(80,-54){$\bullet$}
\put(85,-50){$\overline{1}$} \put(85,-54){$\bullet$}
\put(90,-50){$\overline{1}$} \put(90,-54){$\bullet$}
\put(95,-50){$\overline{1}$} \put(95,-54){$\bullet$}
\put(100,-50){2}             \put(100,-54){$\bullet$}
\put(105,-50){$\overline{2}$} \put(105,-54){$\bullet$}
\put(110,-50){1}              \put(110,-54){$\bullet$}
\put(115,-50){$\overline{1}$}  \put(115,-54){$\bullet$}
\put(25,-60){2} \put(25,-64){$\bullet$}
\put(30,-60){2} \put(30,-64){$\bullet$}
\put(35,-60){2} \put(35,-64){$\bullet$}
\put(40,-60){2}
\put(45,-60){$\overline{2}$} \put(45,-64){$\bullet$}
\put(50,-60){$\overline{2}$} \put(50,-64){$\bullet$}
\put(55,-60){$\overline{2}$} \put(55,-64){$\bullet$}
\put(60,-60){$\overline{2}$}
\put(65,-60){1} \put(65,-64){$\bullet$}
\put(70,-60){1} \put(70,-64){$\bullet$}
\put(75,-60){1} \put(75,-64){$\bullet$}
\put(80,-60){1}
\put(85,-60){$\overline{1}$} \put(85,-64){$\bullet$}
\put(90,-60){$\overline{1}$} \put(90,-64){$\bullet$}
\put(95,-60){$\overline{1}$} \put(95,-64){$\bullet$}
\put(100,-60){$\overline{1}$}
\put(105,-60){2}
\put(110,-60){$\overline{2}$}
\put(115,-60){1}
\put(120,-60){$\overline{1}$}
\put(30,-70){2} \put(30,-74){$\bullet$}
\put(35,-70){2} \put(35,-74){$\bullet$}
\put(40,-70){2}
\put(45,-70){2}
\put(50,-70){$\overline{2}$} \put(50,-74){$\bullet$}
\put(55,-70){$\overline{2}$} \put(55,-74){$\bullet$}
\put(60,-70){$\overline{2}$}
\put(65,-70){$\overline{2}$}
\put(70,-70){1}              \put(70,-74){$\bullet$}
\put(75,-70){1}              \put(75,-74){$\bullet$}
\put(80,-70){1}
\put(85,-70){1}
\put(90,-70){$\overline{1}$} \put(90,-74){$\bullet$}
\put(95,-70){$\overline{1}$} \put(95,-74){$\bullet$}
\put(100,-70){$\overline{1}$}
\put(105,-70){$\overline{1}$}
\put(110,-70){2}
\put(115,-70){$\overline{2}$}
\put(120,-70){1}
\put(125,-70){$\overline{1}$}
\put(35,-80){2} \put(35,-84){$\bullet$}
\put(40,-80){2}
\put(45,-80){2}
\put(50,-80){2}
\put(55,-80){$\overline{2}$} \put(55,-84){$\bullet$}
\put(60,-80){$\overline{2}$}
\put(65,-80){$\overline{2}$}
\put(70,-80){$\overline{2}$}
\put(75,-80){1} \put(75,-84){$\bullet$}
\put(80,-80){1}
\put(85,-80){1}
\put(90,-80){1}
\put(95,-80){$\overline{1}$} \put(95,-84){$\bullet$}
\put(100,-80){$\overline{1}$}
\put(105,-80){$\overline{1}$}
\put(110,-80){$\overline{1}$}
\put(115,-80){2}
\put(120,-80){$\overline{2}$}
\put(125,-80){1}
\put(130,-80){$\overline{1}$} 
\put(40,-90){2}
\put(45,-90){2}
\put(50,-90){2}
\put(55,-90){2}
\put(60,-90){$\overline{2}$}
\put(65,-90){$\overline{2}$}
\put(70,-90){$\overline{2}$}
\put(75,-90){$\overline{2}$}
\put(80,-90){1}
\put(85,-90){1}
\put(90,-90){1}
\put(95,-90){1}
\put(100,-90){$\overline{1}$}
\put(105,-90){$\overline{1}$}
\put(110,-90){$\overline{1}$}
\put(115,-90){$\overline{1}$} \put(115,-94){$\bullet$}
\put(120,-90){2}              \put(120,-94){$\bullet$}
\put(125,-90){$\overline{2}$}
\put(130,-90){1}
\put(135,-90){$\overline{1}$}
\put(45,-100){2}
\put(50,-100){2}
\put(55,-100){2}
\put(60,-100){2}
\put(65,-100){$\overline{2}$}
\put(70,-100){$\overline{2}$}
\put(75,-100){$\overline{2}$}
\put(80,-100){$\overline{2}$}
\put(85,-100){1}
\put(90,-100){1}
\put(95,-100){1}
\put(100,-100){1}
\put(105,-100){$\overline{1}$}
\put(110,-100){$\overline{1}$}
\put(115,-100){$\overline{1}$}  \put(115,-104){$\bullet$}
\put(120,-100){$\overline{1}$}
\put(125,-100){2}                \put(125,-104){$\bullet$}
\put(130,-100){$\overline{2}$}
\put(135,-100){1}
\put(140,-100){$\overline{1}$} 
\put(50,-110){2}
\put(55,-110){2}
\put(60,-110){2}
\put(65,-110){2}
\put(70,-110){$\overline{2}$}
\put(75,-110){$\overline{2}$}
\put(80,-110){$\overline{2}$}
\put(85,-110){$\overline{2}$}
\put(90,-110){1}
\put(95,-110){1}
\put(100,-110){1}
\put(105,-110){1}
\put(110,-110){$\overline{1}$}
\put(115,-110){$\overline{1}$}   \put(115,-114){$\bullet$}
\put(120,-110){$\overline{1}$}
\put(125,-110){$\overline{1}$}
\put(130,-110){2}                \put(130,-114){$\bullet$}
\put(135,-110){$\overline{2}$}
\put(140,-110){1}
\put(145,-110){$\overline{1}$}
\put(55,-120){2}
\put(60,-120){2}
\put(65,-120){2}
\put(70,-120){2}
\put(75,-120){$\overline{2}$}
\put(80,-120){$\overline{2}$}
\put(85,-120){$\overline{2}$}
\put(90,-120){$\overline{2}$}
\put(95,-120){1}
\put(100,-120){1}
\put(105,-120){1}
\put(110,-120){1}               \put(110,-124){$\bullet$}
\put(115,-120){$\overline{1}$}  \put(115,-124){$\bullet$}
\put(120,-120){$\overline{1}$}
\put(125,-120){$\overline{1}$}
\put(130,-120){$\overline{1}$}
\put(135,-120){2}               \put(135,-124){$\bullet$}
\put(140,-120){$\overline{2}$}  \put(140,-124){$\bullet$}
\put(145,-120){1}
\put(150,-120){$\overline{1}$} 
\put(60,-130){2}
\put(65,-130){2}
\put(70,-130){2}
\put(75,-130){2}
\put(80,-130){$\overline{2}$}
\put(85,-130){$\overline{2}$}
\put(90,-130){$\overline{2}$}
\put(95,-130){$\overline{2}$}
\put(100,-130){1}
\put(105,-130){1}
\put(110,-130){1}                \put(110,-134){$\bullet$}
\put(115,-130){1}
\put(120,-130){$\overline{1}$}
\put(125,-130){$\overline{1}$}
\put(130,-130){$\overline{1}$}
\put(135,-130){$\overline{1}$}
\put(140,-130){2}
\put(145,-130){$\overline{2}$}  \put(145,-134){$\bullet$}
\put(150,-130){1} 
\put(155,-130){$\overline{1}$}
\put(65,-140){2}
\put(70,-140){2}
\put(75,-140){2}
\put(80,-140){2}
\put(85,-140){$\overline{2}$}
\put(90,-140){$\overline{2}$}
\put(95,-140){$\overline{2}$}
\put(100,-140){$\overline{2}$}
\put(105,-140){1}
\put(110,-140){1}              \put(110,-144){$\bullet$}
\put(115,-140){1}
\put(120,-140){1}
\put(125,-140){$\overline{1}$}
\put(130,-140){$\overline{1}$}
\put(135,-140){$\overline{1}$}
\put(140,-140){$\overline{1}$}
\put(145,-140){2}
\put(150,-140){$\overline{2}$}  \put(150,-144){$\bullet$}
\put(155,-140){1}
\put(160,-140){$\overline{1}$} 
\put(70,-150){2}
\put(75,-150){2}
\put(80,-150){2}
\put(85,-150){2}
\put(90,-150){$\overline{2}$}
\put(95,-150){$\overline{2}$}
\put(100,-150){$\overline{2}$}
\put(105,-150){$\overline{2}$}  \put(105,-154){$\bullet$}
\put(110,-150){1}                \put(110,-154){$\bullet$}
\put(115,-150){1}
\put(120,-150){1}
\put(125,-150){1}
\put(130,-150){$\overline{1}$}
\put(135,-150){$\overline{1}$}
\put(140,-150){$\overline{1}$}
\put(145,-150){$\overline{1}$}
\put(150,-150){2}
\put(155,-150){$\overline{2}$}  \put(155,-154){$\bullet$}
\put(160,-150){1}               \put(160,-154){$\bullet$}
\put(165,-150){$\overline{1}$}
\put(75,-160){2}
\put(80,-160){2}
\put(85,-160){2}
\put(90,-160){2}
\put(95,-160){$\overline{2}$}
\put(100,-160){$\overline{2}$}
\put(105,-160){$\overline{2}$} \put(105,-164){$\bullet$}
\put(110,-160){$\overline{2}$}
\put(115,-160){1}
\put(120,-160){1}
\put(125,-160){1}
\put(130,-160){1}
\put(135,-160){$\overline{1}$}
\put(140,-160){$\overline{1}$}
\put(145,-160){$\overline{1}$}
\put(150,-160){$\overline{1}$}
\put(155,-160){2}
\put(160,-160){$\overline{2}$}
\put(165,-160){1}             \put(165,-164){$\bullet$}
\put(170,-160){$\overline{1}$} 
\put(80,-170){2}
\put(85,-170){2}
\put(90,-170){2}
\put(95,-170){2}
\put(100,-170){$\overline{2}$}
\put(105,-170){$\overline{2}$}  \put(105,-174){$\bullet$}
\put(110,-170){$\overline{2}$}
\put(115,-170){$\overline{2}$}
\put(120,-170){1}
\put(125,-170){1}
\put(130,-170){1}
\put(135,-170){1}
\put(140,-170){$\overline{1}$}
\put(145,-170){$\overline{1}$}
\put(150,-170){$\overline{1}$}
\put(155,-170){$\overline{1}$}
\put(160,-170){2}
\put(165,-170){$\overline{2}$}
\put(170,-170){1}               \put(170,-174){$\bullet$}
\put(175,-170){$\overline{1}$}
\put(85,-180){2}
\put(90,-180){2}
\put(95,-180){2}
\put(100,-180){2}                \put(100,-184){$\bullet$}
\put(105,-180){$\overline{2}$}     \put(105,-184){$\bullet$}
\put(110,-180){$\overline{2}$}
\put(115,-180){$\overline{2}$}
\put(120,-180){$\overline{2}$}
\put(125,-180){1}
\put(130,-180){1}
\put(135,-180){1}
\put(140,-180){1}
\put(145,-180){$\overline{1}$}
\put(150,-180){$\overline{1}$}
\put(155,-180){$\overline{1}$}
\put(160,-180){$\overline{1}$}
\put(165,-180){2}
\put(170,-180){$\overline{2}$}
\put(175,-180){1}                \put(175,-184){$\bullet$}
\put(180,-180){$\overline{1}$}   \put(180,-184){$\bullet$}
\put(90,-190){2}
\put(95,-190){2}
\put(100,-190){2}                  \put(100,-194){$\bullet$}
\put(105,-190){2}
\put(110,-190){$\overline{2}$}
\put(115,-190){$\overline{2}$}
\put(120,-190){$\overline{2}$}
\put(125,-190){$\overline{2}$}
\put(130,-190){1}
\put(135,-190){1}
\put(140,-190){1}
\put(145,-190){1}
\put(150,-190){$\overline{1}$}
\put(155,-190){$\overline{1}$}
\put(160,-190){$\overline{1}$}
\put(165,-190){$\overline{1}$}
\put(170,-190){2}
\put(175,-190){$\overline{2}$}
\put(180,-190){1}
\put(185,-190){$\overline{1}$}     \put(185,-194){$\bullet$}
\put(95,-200){2}
\put(100,-200){2}                   \put(100,-204){$\bullet$}
\put(105,-200){2}
\put(110,-200){2}
\put(115,-200){$\overline{2}$}
\put(120,-200){$\overline{2}$}
\put(125,-200){$\overline{2}$}
\put(130,-200){$\overline{2}$}
\put(135,-200){1}
\put(140,-200){1}
\put(145,-200){1}
\put(150,-200){1}
\put(155,-200){$\overline{1}$}
\put(160,-200){$\overline{1}$}
\put(165,-200){$\overline{1}$}
\put(170,-200){$\overline{1}$}
\put(175,-200){2}
\put(180,-200){$\overline{2}$}
\put(185,-200){1}
\put(190,-200){$\overline{1}$}    \put(190,-204){$\bullet$}
\end{picture}

\caption{$W_0^2$ cannot cover $W_0^{-1}$. The first row is $W_0^2$.
The second to the last rows describe all the possible positions for 
$W_0^{-1}$.
We put $\bullet$ below illegal pairs of $(1,2)$ or 
$(\overline{1}, \overline{2})$.
In each position for $W_0^{-1}$, we find at least one such pair, 
which is a contradiction.}
\end{center}
\end{figure}

\noindent
(3) Observe that $w_j(1,+)$ consists of $a_1$ and $b$,
while $w_i$ contains $a_2$.
If $w_j(1,+)$ covers $w_i$, then each $a_2$ in $w_i$ faces some 
$a_1$ in $w_j(1,+)$, which contradicts (1).
The same argument goes for $w_j(1,-)^{-1}$, and we omit it.

\noindent
(4) Assume $w_i^k$ covers $w_j$ for some $i <j$ and some $k > 0$.
Then $w_j(1,+)$ covers some $w_i$ of $w_i^k$ since 
$2|w_i| \le |w_j(1,+)|$, which contradicts (3).
Assume $w_i^k$ covers $w_j^{-1}$ for some $i<j$ and some $k > 0$.
Then $w_j(1,-)^{-1}$ covers some $w_i$ of $w_i^k$,
which contradicts (3).
\qed

\begin{Lemma}
$c_{w_i}$ and $c_{w_i^{-1}}$ satisfy the following conditions.

\noindent
(1) For any $n \ge 1$ and any $i \ge 0$, 
$c_{w_i} (w_i^n)=n.$

\noindent
(2) For any $n \ge 1$ and any $i \ge 0$, 
$c_{w_i^{-1}}(w_i^n) =0.$

\noindent
(3) For any $j > i \ge0$ and any $n \ge 1$, 
$c_{w_j}(w_i^n)=0.$

\noindent 
(4) For any $j > i \ge0$ and any $n \ge 1$,
$c_{w_j^{-1}}(w_i^n) =0.$
\end{Lemma}

\noindent
{\it Proof}. (1) Since $w_i^n$ is reduced, $w_i^n$ is a geodesic.
Thus $c_{w_i}(w_i^n) \ge |w_i^n|_{w_i} =n$.
On the other hand, $c_{w_i}(w_i^n) \le |\overline{w_i^n}|/|w_i| =n$.
Thus $c_{w_i}(w_i^n) =n.$

\noindent
(2) To show  $c_{w_i^{-1}}(w_i^n) =0$ by contradiction, 
assume $c_{w_i^{-1}}(w_i^n) > 0$. Take a realizing geodesic of 
$c_{w_i^{-1}}$ at $w_i^n$, $\alpha$. Then $|\alpha|_{w_i^{-1}} > 0$. 
Fix a subword labeled by $w_i^{-1}$ in $\alpha$. 
Since $w_i^n$ and $\alpha$ are reduced, there is a subword of $w_i^n$ 
labeled by $w_i^2$ which covers the $w_i^{-1}$ in $\alpha$. 
But $w_i^2$ cannot cover $w_i^{-1}$ by Lemma 4.2(2), a contradiction.
See Figure 7.

\begin{figure}
\centerline{
\epsfbox{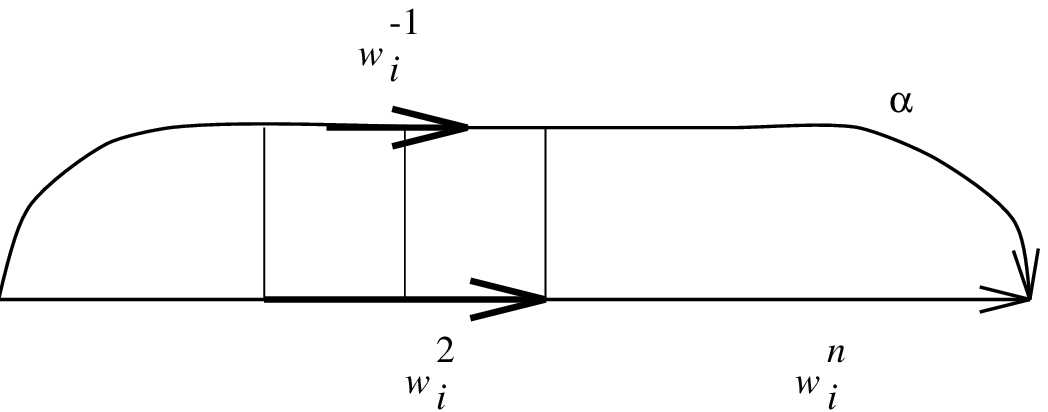}
}
\caption{$w_i^2$ cannot cover $w_i^{-1}$.}
\end{figure}

\noindent
(3) To show the claim by contradiction, 
assume $c_{w_j}(w_i^n) > 0$ for some $i < j$ and $n$. 
Take a realizing geodesic $\alpha$ of $c_{w_j}$ at $w_i^n$. 
Then $|\alpha|_{w_j} >0$. Fix some subword labeled by $w_j$ in $\alpha$. 

Since $\alpha$ and $w_i^n$ are reduced, there exists a subword labeled by 
$w_i^k$ in $w_i^n$ which covers the $w_j$ in $\alpha$. 
But $w_i^k$ cannot cover $w_j$ by Lemma 4.2(4), a contradiction.

\noindent
(4) Similar to (3).
\qed

Now we are in a position to prove Proposition 2.

\noindent
{\it Proof of Proposition 2}. 
Let $w_i, 0 \le i < \infty$ be the set of words defined before in this section.
We will show they are the words we want.
Since $h_{w_i} = c_{w_i} -c_{w_i^{-1}}$, 
$h_{w_i}(w_i^n)=n$ by Lemma 4.3(1) and (2). 
By Lemma 4.3(3) and (4), $h_{w_j}(w_i^n)=0$ for $i<j$. 
For a homomorphism $\phi:G \to \mathbf{R} $, $\phi(\overline{w_i}) =0$ since 
$\overline{w_i} \in [G,G]$ by Lemma 4.1(3).
\qed

\section{Proofs of Th 1 and Cor 1,2,3}

\noindent
{\it Proof of Th 1}. 
Let $w_i, \, 0 \le i < \infty$ be the words in Prop 2.
By Prop 1, all the 
cocycles $\delta h_{w_i}$ have the same bound. 
It follows that if $(a_i)_i \in l^1$, then 
$\sum_i a_i \delta h_{w_i}$ is also a well-defined cocycle.
We get a real linear map
$$\omega: l^1 \to H^2_b(G;\mathbf R), $$
which sends $(a_i)_i$ to the cohomology class of $\sum_i a_i \delta h_{w_i}$. 
In order to show $\omega$ is injective, suppose $\omega ((a_i))=0$.  
Then 
$$\delta(\sum_{i=0}^{\infty} a_i h_{w_i}) =\delta b$$
for some $b \in C_b^1(G;\mathbf R)$. This means 
$$\sum _i a_i h_{w_i} -b = \phi,$$
for some homomorphism $\phi:G \to \mathbf R$. 
Applying this to $\overline{w_0^n} \in G$, we find 
$$a_0 n - b(\overline{w_0^n}) = \phi(\overline{w_0^n}) =0,$$
for all $n \ge 1$ by Prop 2. 
Since $b$ is bounded, $a_0=0$. 
Similarly, $a_i=0$ for all $i \ge 1$. Thus $\omega$ is injective.
It is well-known that the cardinality of the dimension of $l^1$
as a vector space is continuum. \qed

\noindent
{\it Proof of Cor 1}. 
W.l.o.g., $|A| \ge 3$.
Since $C=\{ 1 \}$, $|C \backslash A / C| = |A| \ge 3$.
Apply Th 1. \qed

\noindent
{\it Proof of Cor 2}.
Since $|C|  < \infty$ and $|A|=\infty$, $|C \backslash A / C| = \infty$.
Apply Th 1.
\qed

\noindent
{\it Proof of Cor 3}.
Since $A$ is abelian, 
$|C \backslash A / C| = |A/C| \ge 3$.
Apply Th 1. \qed

\section{Quasi homomorphisms of $A*_C \varphi$}

Let $G=A*_C \varphi = < A,t; c=t^{-1}\varphi(c)t, \, ^{\forall} c \in C>$
with $|A/C| \ge 2, \, |A/ \varphi (C)| \ge 2$. 
Let $g$ be 
$$g =a_1 t^{n_1} a_2 t^{n_2} \dots a_I t^{n_I} a_{I+1}, $$
with $a_i \in A$ and $ n_i \neq 0$, where $a_1, a_{I+1}$ may be empty.
Suppose $g$ satisfies the following two conditions (i) and (ii)
for $1 \le ^{\forall} i \le I-1$.

\noindent
(i) If $n_i >0, n_{i+1} <0$, then $a_{i+1}\notin C$.

\noindent
(ii) If $n_i < 0, n_{i+1} > 0$, then $a_{i+1} \notin \varphi(C)$.

\noindent
Then we say $g$ is {\it reduced}.

The following fact is known as Britton's lemma.

\begin{Lemma} 
[{\bf Britton, \cite{LS}}] 
Suppose $1 \le I$. If $g$ is reduced, then $g \neq 1$ in $G$.
\end{Lemma}

As an application of Britton's lemma, we have:

\begin{Lemma}
Let 
$$g=a_1 t^{n_1} a_2 t^{n_2} \dots a_I t^{n_I} a_{I+1}, \, \, a_i \in A$$
$$h=b_1 t^{m_1} b_2 t^{m_2} \dots b_J t^{m_J} b_{J+1}, \, \, b_j \in A$$
be reduced with $n_i = \pm 1, \, m_j = \pm 1.$
If $g=h$ in $G$, then
$$I=J, \, n_1=m_1, \dots ,n_I=m_I,$$
and 
$$a_{I+1}b_{I+1}^{-1}, \, \,  t^{n_I}a_{I+1}b_{I+1}^{-1}t^{-m_I}, \, \, 
a_It^{n_I}a_{I+1}b_{I+1}^{-1}t^{-m_I}b_I^{-1}, \dots $$
$$\dots , t^{n_1}a_2t^{n_2}a_3  \dots a_{I+1}b_{I+1}^{-1} \dots 
b_3 ^{-1}t^{-m_2}b_2^{-1}t^{-m_1}$$
are in $C$ or $\varphi(C)$, especially, in $A$.
\end{Lemma}

\noindent
{\it Proof}. W.l.o.g., we assume $J \le I$.
$h^{-1}= b_{J+1}^{-1}t^{-m_J}b^{-1}_J t^{-m_{J-1}} \dots b^{-1}_1$
is reduced since $h$ is reduced.
Put 
$$d_{I+1} =a_{I+1} b_{J+1}^{-1}, $$
then $d_{I+1} \in A$ and 

$$gh^{-1} = (a_1 t^{n_1} a_2 t^{n_2} \dots a_I t^{n_I}) d_{I+1} 
            (t^{-m_J}b^{-1}_J t^{-m_{J-1}} \dots b^{-1}_1)$$
in $G$. 
Since $gh^{-1} =1$ in $G$, the expression on the right hand side 
is not reduced by Britton's lemma.
But, since $a_1 t^{n_1} a_2 t^{n_2} \dots a_I t^{n_I}$ and 
$t^{-m_J}b^{-1}_J t^{-m_{J-1}} \dots b^{-1}_1$ are reduced, we have 
$$n_I =m_J=1 (\mbox{or} -1), \, \, d_{I+1} \in C (\varphi(C), \mbox{resp.}),$$
and $t^{n_I} d_{I+1} t^{-m_J} \in \varphi(C)$($C$, resp.).

Put
$$c_I = t^{n_I} d_{I+1} t^{-m_J}, \, \, d_I = a_I c_I b_J^{-1}.$$
Clearly $d_I \in A$ and 
$$gh^{-1} =  (a_1 t^{n_1} a_2 t^{n_2} \dots a_{I-1} t^{n_{I-1}}) d_I 
            (t^{-m_{J-1}}b^{-1}_{J-1} t^{-m_{J-2}} \dots b^{-1}_1)$$
in $G$. 
By Britton's lemma, the expression on the right hand side is not reduced. 
Thus $$n_{I-1}=m_{J-1}=1(\mbox{or} -1), $$
$$d_I \in C (\varphi(C), \mbox{resp}.), \, \,  
\mbox{and } c_{I-1}=t^{n_{I-1}} d_I t^{-m_{J-1}} 
\in \varphi(C) (C, \mbox{resp}.).$$

We define $d_{I-1}, d_{I-2}, \dots, d_{I-J+2}$ and 
$c_{I-1}, c_{I-2}, \dots, c_{I-J+1}$ inductively by 
$$c_i=t^{n_i} d_{i+1} t^{-m_{J-I+i}}, \, \, \, 
  d_i = a_i c_i b^{-1}_{J-I+i}.$$
See Figure 8. 

\begin{figure}
\centerline{\epsfxsize=12cm
\epsfbox{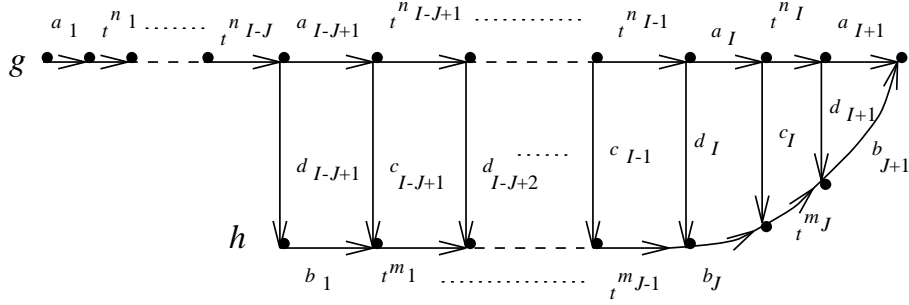}
}
\caption{This illustrates Lemma 6.2. 
$d_i, c_i \in A$ for $I-J+1 \le ^{\forall} i \le I$.  
If $J < I$, then it contradicts Britton's lemma.}
\end{figure}

Using Britton's lemma repeatedly, we have 
$$n_{I-J+j}=m_j, \, 1 \le ^{\forall} j \le J$$ 
 and
$$d_{I+1}, d_I, \dots, d_{I-J+2}, c_I, c_{I-1}, \dots, c_{I-J+1} 
\in C \mbox{ or } \varphi(C).$$
Put 
$$d_{I-J+1} = a_{I-J+1} c_{I-J+1} b_1^{-1}.$$
Clearly $d_{I-J+1} \in A$.
To complete the proof, it suffices to show $J=I$.
In order to show it by contradiction, suppose $J < I$. Then 
$$gh^{-1} = a_1t^{n_1}a_2 \dots t^{n_{I-J-1}} a_{I-J}t^{n_{I-J}} d_{I-J+1} $$
in $G$ and the expression on the right hand side is reduced, 
which contradicts Britton's 
lemma, since $gh^{-1}=1$ in $G$. We got $I=J$.
\qed

We take $\{ t\} \cup  A \backslash \{ 1 \}$ as a set of generators 
of $G$ and write the Cayley graph of $G$ w.r.t. this set by $\Gamma$.

\begin{Lemma}
If a path $\alpha$ is a geodesic in $\Gamma$, then it is reduced.
\end{Lemma}

\noindent
{\it Proof}.
If a path $\alpha$ is not reduced, then we can make it shorter using 
$tct^{-1} = \varphi(c)$, hance $\alpha$ is not a geodesic.
\qed

\noindent
{\it Remark}. Compare Lemma 6.3 with Lemma 3.1. 
A reduced path is not always a geodesic in this case.
For example, the left hand side of $\varphi(c)^{-1} t c =t$ for 	
$c \in C \backslash \{1\}$ is reduced but not a geodesic.

\begin{Lemma}
Let $\alpha$ be a path and $w$ a word. If $w^2$ is reduced, 
then there is a reduced path $\beta$ which realizes $c_w$ at $\alpha$.
\end{Lemma}

\noindent
{\it Proof}. Take a path $\beta$ which realizes $c_w$ at $\alpha$ s.t.
$|\beta|_w$ is minimal among all the realizing paths. Then a similar 
argument to the proof of Lemma 3.2 shows $\beta$ is reduced. 
We omit the detail.
\qed

\begin{Lemma}
Let $\alpha, \beta$ be paths starting at $1$. We have 
$$|c_w(\alpha)-c_w(\beta)| \le 2|\overline{\alpha^{-1} \beta}|.$$
$$|h_w(\alpha)-h_w(\beta)| \le 4|\overline{\alpha^{-1} \beta}|.$$

\end{Lemma}

\noindent
{\it Proof}. Similar to Lemma 3.3 and 3.4. \qed

\begin{Lemma}
Let $\alpha$ be a reduced path. If $\alpha= \alpha_1 \alpha_2$, then 
$$|h_w(\alpha) - h_w(\alpha_1) -h_w(\alpha_2) | \le 10.$$
\end{Lemma}

\noindent
{\it Proof}. Similar to Lemma 3.7. \qed

\begin{Prop}
$$|\delta h_w | \le 78.$$
\end{Prop}

\noindent
{\it Proof}. The outline is similar to the proof of Proposition 1.
Let $x,y$ be elements in $G$. We show $|h_w(xy)-h_w(x)-h_w(y)| \le 78.$
Take reduced paths $\alpha, \beta, \gamma$ with 
$\overline{\alpha}=x, \overline{\beta}=y, \overline{\gamma}=xy$.
Since $\overline{\alpha \beta} = \overline{\gamma}$ and 
$\alpha, \beta, \gamma$ are 
reduced, by Lemma 6.1 and 6.2, there exist subdivisions of 
$\alpha, \beta, \gamma$,
$$\alpha=\alpha_1 \alpha_2 \alpha_3, \, 
\beta =\beta_1 \beta_2 \beta_3, \, \gamma=\gamma_1 \gamma_2 \gamma_3, $$
s.t. $\overline{\gamma^{-1}_1 \alpha_1}= c_1, 
     \overline{\alpha_3 \beta_1}=c_2, \overline{\beta_3 \gamma_3^{-1}}=c_3$, 
for some $c_1, c_2, c_3 \in C \cup \varphi(C),$ and 
$\overline{\alpha_2}, \overline{\beta_2}, \overline{\gamma_2} \in A$. 
By Lemma 6.6, 
$$|h_w(\alpha) -h_w(\alpha_1) -h_w(\alpha_2) -h_w(\alpha_3) | \le 20.$$
Since $|h_w(\alpha_2) | \le 2,$
$|h_w(\alpha) -h_w(\alpha_1) - h_w(\alpha_3) | \le 22.$
Similarly, $$|h_w(\beta) -h_w(\beta_1) - h_w(\beta_3) | \le 22, \, 
|h_w(\gamma) -h_w(\gamma_1) - h_w(\gamma_3) | \le 22.$$
By a similar argument to the proof of Prop 1, 
$$|h_w(x)+h_w(y)-h_w(xy)| = |h_w(\alpha) + h_w(\beta) -h_w(\gamma)| \le 78.$$
\qed

\section{Choice of words for $A*_C\varphi$}
Let $G=A*_C \varphi$ with $|A/C| \ge 2, \, |A/ \varphi (C)| \ge 2$.

\begin{Lemma}
Let $w$ be 
$$w=t^{n_1} a_1 t^{n_2} a_2 \dots t^{n_I} a_I,$$
s.t. $a_i \in A \backslash \{ 1 \}, n_i  \in {\mathbf Z} \backslash \{0\}$
for $1 \le ^{\forall} i \le I$. 
We denote the set of the following conditions (1.1), \dots, (1.4) 
by Condition I and (2.1), \dots, (2.4) by Condition II.

\noindent
$(1.1) \,  0 < n_1, n_3, n_5, \dots $.

\noindent
$(1.2) \,0 > n_2, n_4, n_6, \dots$.

\noindent
$(1.3) \,  a_1, a_3, a_5, \dots \notin C$.

\noindent
$(1.4) \, a_2, a_4, a_6, \dots \notin \varphi(C)$.

\noindent
$(2.1) \, 0 > n_1, n_3, n_5, \dots $.

\noindent
$(2.2) \, 0 < n_2, n_4, n_6, \dots$. 

\noindent
$(2.3) \, a_1, a_3, a_5, \dots \notin \varphi(C)$.

\noindent
$(2.4) \,a_2, a_4, a_6, \dots \notin C$.

If either the Condition I or II holds, 
then $w$ is a geodesic in $\Gamma$.
\end{Lemma}

\noindent
{\it Proof}.
We treat the case that the Condition I holds.
The other case is similar.
Under this condition, clearly $w$ is reduced.
Let $\gamma$ be a geodesic with $\overline{w} = \overline{\gamma}$.
Then by Lemma 6.3, $\gamma$ is reduced. By Lemma 6.2, we have 
$$\gamma =b_0 \tau_1 b_1 \tau_2 b_2 
\dots \tau_I b_I,$$ 
s.t. $b_i \in A \backslash \{ 1 \}$ for $1 \le i \le I-1$, 
and $b_0, b_I \in  A \backslash \{ 1 \}$ or is empty, and 

\begin{eqnarray*}
\tau_i & = & t b_{i,1} t b_{i,2} \dots t b_{i,n_i -2} t b_{i,n_i-1} t,  
\mbox{   if $i$ is odd}, \\
\tau_i & = & t^{-1} b_{i,1} t^{-1} b_{i,2} \dots t^{-1} b_{i,n_i -2} 
t^{-1} b_{i,n_i-1} t^{-1},  
\mbox{   if $i$ is even}, 
\end{eqnarray*}
where $b_{i,j} \in A \backslash \{ 1 \}$ or is empty, 
for $1 \le ^{\forall} i \le I$ and $1 \le ^{\forall} j \le n_i -1$. 

We will show $b_I$ is not empty. To show this by contradiction, 
assume $b_I$ is empty.
Then

\begin{eqnarray*}
w \gamma^{-1} & = & \dots t a_I t^{-1} \dots, \mbox{   if $I$ is odd}, \\
w \gamma^{-1}    & = & \dots t^{-1} a_I t \dots, \mbox{   if $I$ is even.}
\end{eqnarray*}
Since $w$ and $\gamma^{-1}$ are reduced, $w \gamma^{-1}$ is reduced
by (1.3) if $I$ is odd or by (1.4) if $I$ is even. 
Then it follows from Britton's lemma that 
$\overline{w \gamma^{-1}} \neq 1$, a contradiction.
We got $b_I$ is not empty. 
Thus 
$$|\gamma| \ge \sum_{i=1}^I |\tau_i| + I \ge \sum_{i=1}^I |n_i| +I = |w|.$$
We have the first inequality because $b_1, \dots, b_I$ are not empty, and 
the second one since the number of $t$'s in $\tau_i$ is $n_i$ 
for $1 \le ^{\forall} i \le I$.
Since $\gamma$ is a geodesic, we have $|\gamma| = |w|$, 
hence $w$ is a geodesic.
\qed

\begin{Lemma}
Take $g \in A \backslash C$ and $h \in A \backslash \varphi(C)$ and fix them.

\noindent
Let $w_i, 0 \le i < \infty$, be words s.t.
$$w_i = t^{10^i} g t^{-10^i} h 
        t^{10^i} g^{-1} t^{-10^i} h^{-1}
        t^{2 \cdot 10^i} g t^{-2 \cdot 10^i} h 
        t^{3 \cdot 10^i} g^{-1} t^{-3 \cdot 10^i} h^{-1}.$$
Then the words satisfy the following conditions. 

\noindent
(1) For any $i \ge 0$ and any $n \ge 1$, 
$c_{w_i} (w_i^n) =n.$

\noindent
(2) For any  $i \ge 0$ and any $n \ge 1$,
$c_{w_i^{-1}} (w_i^n) =0.$

\noindent
(3) For any  $j > i \ge 0$ and any $n \ge 1$,
$c_{w_j^{\pm 1}} (w_i^n) =0.$

\noindent
(4)  For any  $i \ge 0$, 
$\overline{w_i} \in [G,G]$.
\end{Lemma}

\noindent
{\it Proof}.
(1) By Lemma 7.1, $w_i^n$ is a geodesic. Thus $c_{w_i} (w_i^n) =n$.

\noindent
(2) In order to show $c_{w_0^{-1}}(w_0^n) =0$ for any $n \ge 1$ 
by contradiction, 
suppose $c_{w_0^{-1}}(w_0^n) > 0$ for some $n$.
By Lemma 6.4, take a reduced path $\alpha$ s.t. $\overline{\alpha} =
\overline{w_0^n}$
and $|\alpha|_{w_0^{-1}} > 0$. 
We describe the order of $t$'s and $t^{-1}$'s appearing in $w_0$ by $W_0$,
$$W_0 = +\, -\, +\, -\, +\, +\, -\, -\, +\, +\, +\, -\, -\, -,$$
where $+$ stands for $t$, $-$ for $t^{-1}$ 
and we ignore $g^{\pm 1}, h^{\pm 1}$ in $w_0$. 
Then $W_0^n$ represents $w_0^n$. 
Let $\Lambda $ be the order of $t$'s and $t^{-1}$'s in $\alpha$. 
Applying Lemma 6.2 for $\alpha$ and $w_0^n$, we get $W_0^n= \Lambda$.
$\Lambda$ contains $W_0^{-1}$ as a subpiece since $|\alpha|_{w_0^{-1}} > 0$,
where 
$$W_0^{-1} = +\, +\, +\, -\, -\, -\, +\, +\, -\, -\, +\, -\, +\, -\, .$$
But it is easy to check that $W_0^n$ cannot contain $W_0^{-1}$ as a subset, 
thus we get a contradiction. Hence $c_{w_0^{-1}}(w_0^n) =0$.
Similarly, we can show that $c_{w_i^{-1}}(w_i^n) =0$ 
for any $i \ge 1$ and any $n \ge 1$.

\noindent 
(3) Let $W_i, \, 1 \le i < \infty$ be the order of $t$'s 
and $t^{-1}$'s in $w_i$.
To show $c_{w_j}(w_i^n) =0$ for all $n \ge 1$ and all $j > i >0$
by contradiction, 
assume $c_{w_j}(w_i^n) > 0$ for some $n \ge 1$ and some $j > i >0$.
Then, by our assumption, $W_j$ must be contained in $W_i^n$ as a subword.
This is a contradiction since $W_j$ has $3 \cdot 10^j$ consecutive $t$'s,
but $W_i$ has at most $3 \cdot 10^i$ consecutive $t$'s.
Thus we have $c_{w_j}(w_i^n) =0$. 
$c_{w_j^{-1}}(w_i^n) =0$ is similar.

\noindent
(4) Set $T=t^{10^i}$, then 
$$w_i=[T,g][gh,[T,g^{-1}]\, ][T,g^{-1}][g,h][T^2,g][gh,[T^3,g^{-1}]\, ]
         [T^3,g^{-1}][g,h],$$
thus $\overline{w_i} \in [G,G]$.
\qed

\begin{Prop}
Let $G=A*_C \varphi$ with $|A/C| \ge 2, \, |A / \varphi(C)| \ge 2$.
There exist words $w_i, \, 0 \le i < \infty$ which satisfy the 
following conditions.

\noindent
(1) For any $i \ge 0$ and any $n \ge 1$, 
$h_{w_i} (w_i^n) =n. $

\noindent
(2) For any  $j > i \ge 0$ and any $n \ge 1$,
$h_{w_j} (w_i^n) =0.$

\noindent
(3) For any $i \ge 0$, 
$\overline{w_i} \in [G,G]$.
\end{Prop}

\noindent
{\it Proof}. The words $w_i, \, 0 \le i < \infty$ in Lemma 7.2
clearly satisfy (1), (2) and (3). \qed

\section{Proofs of Th 2 and Th 3}

\noindent
{\it Proof of Th 2}.
We can show Th 2 using Prop 3 and Prop 4 as we showed 
Th 1 using Prop 1 and Prop 2. We omit the detail. 
\qed

\noindent
{\it Proof of Th 3}.
By J. Stallings' structure theorem \cite{S}, $G$ is either

\noindent
(1) $A*_C B$ with $|C| < \infty, \, |A/C| \ge 3$ and $|B/C| \ge 2$,

\noindent
or

\noindent
(2) $A *_C \varphi$ with $|C| <  \infty$ and $|A/C| \ge 2$.

In the case (1), if $|A| = \infty$ or $|B|=\infty$, then we have the conclusion
by Cor 2. If $|A| < \infty$ and $|B| < \infty$, then $G$ is word-hyperbolic.
Since $G$ has infinitely many ends, it is non-elementary.
It is known that the conclusion of Th 3 holds for a non-elementary 
word-hyperbolic group, \cite{EF}. 
In the case (2), since $|A/C| \ge 2$ and $|C| < \infty$, 
we have $|A/ \varphi (C)| \ge 2$. Apply Th 2. 
\qed


\begin{thebibliography}{BaGh3}

\bibitem[BaGh1]{BaGh1}
J. Barge, E. Ghys, Surfaces et cohomologie born\'{e}e, 
Invent. Math., 92, 1988, 509-526.

\bibitem[BaGh2]{BaGh2}
J. Barge, E. Ghys, Cocycles born\'{e}s et actions de groupes 
sur les arbres r\'{e}els,
in ``Group Theory from a Geometric Viewpoint'', 
World Sci. Pub., New Jersey, 1991, 617-622.

\bibitem[BaGh3]{BaGh3}
J. Barge, E. Ghys, Cocycles d'Euler et de Maslov, Math. Ann., 294, 
no 2, 1992, 235-265.

\bibitem[Bav]{Bav}
C. Bavard, Longueur stable des commutateurs, L'Enseignement Math., 37, 1991, 
 109-150.

\bibitem[B]{B}
R. Brooks, 
Some remarks on bounded cohomology, 
Ann. Math. Studies, 97, 1881, 53-63.

\bibitem[BS]{BS}
R. Brooks, C. Series,
Bounded cohomology for surface groups, Topology, 23, no 1, 1884, 29-36.

\bibitem[Br]{Br}
K.S. Brown, ``Cohomology of groups'', Springer, New York, 1982.


\bibitem[EF]{EF} D.B.A. Epstein, K. Fujiwara, The second bounded 
cohomology of word-hyperbolic groups, preprint.

\bibitem[Gh]{Gh} E. Ghys, Groupes d'hom\'{e}omorphismes du cercle 
et cohomologie born\'{e}e, Contemp. Math. 58, AMS, 1987, 81-106.

\bibitem[Gr]{Gr} R.I. Grigorchuk, Some remarks on bounded cohomology, LMS
Lecture Note 204, Cambridge Univ Press, Cambridge, 1995, 111-163.

\bibitem[G]{G}
M. Gromov, 
Volume and bounded cohomology, Publ. Math. IHES, 56, 1982, 5-99.

\bibitem[I]{I}
N.V. Ivanov,
Foundations of the theory of bounded cohomology, 
Zap. Nauchn. Sem. Leningr. Otd. Mat. Inst., 143, 1985, 69-109.

\bibitem[LS]{LS}
R.C. Lyndon, P.E. Schupp, ``Combinatorial Group Theory'', Berlin, 
Springer, 1977.

\bibitem[MaMo]{MaMo}
S. Matsumoto, S. Morita, Bounded cohomology of certain groups 
              of homomorphisms, Proc. A.M.S., 94, 1985, 539-544.

\bibitem[Mi]{Mi}
Y. Mitsumatsu, Bounded cohomology and $l^1$-homology of surfaces, 
Topology, 23, 1984, 465-471.



\bibitem[So1]{So1}
T. Soma, Bounded cohomology of closed surfaces, preprint.

\bibitem[So2]{So2}
T. Soma, Bounded cohomology and topologically tame Kleinian groups, preprint.


\bibitem[S]{S}
J. Stallings, ``Group theory and three-dimensional manifolds'', 
                 Yale Univ. Press, New Haven, 1971.

\bibitem[Y]{Y}
T. Yoshida, On 3-dimensional bounded cohomology of surfaces, in ``Homotopy 
Theory and Related Topics'', Advanced Studies in Pure Math 9, Kinokuniya,
Tokyo, 1986, 173-176.

\end{thebibliography}
\end{document}